\documentclass[11pt,leqno]{article}

\usepackage{amssymb,amsmath,amsthm}
\newcommand{\Int}{{\rm Int}}
\newcommand{\Ind}{{\rm Ind}}

\newcommand{\rmpt}{{\rm pt}}
\newcommand{\Set}{{\bf Set}}
\newcommand{\supp}{{\rm supp}}

\newcommand{\shcw}{\shc^\wedge}

\newcommand{\eq}{\begin{eqnarray}}
\newcommand{\eneq}{\end{eqnarray}}
\newcommand{\eqn}{\begin{eqnarray*}}
\newcommand{\eneqn}{\end{eqnarray*}}

\newcommand{\indlim}[1][]{\mathop{\varinjlim}\limits_{#1}}
\newcommand{\prolim}[1][]{\mathop{\varprojlim}\limits_{#1}}

\newcommand{\indc}{{\rm Ind}(\shc)}

\newcommand{\inddlim}[1][]{\mathop{``{\varinjlim}"}\limits_{#1}}

\newcommand{\llim}[1][]{\mathop{lim\,}\limits_{#1}}

\def\etens{\mathbin{\boxtimes}}
\newcount\temp
\temp=\catcode`\@
\catcode`\@=11
\def\@bletens{\mathbin{\etens^{L}}}
\def\@letens_#1{\mathbin{\etens_{\raise1.5ex\hbox to-.1em{}#1}^{L}}}
\def\letens{\@ifnextchar _{\@letens}{\@bletens}}
\catcode`\@=\temp

\def\rop{{\rm op}}

\def\phi{{\varphi}}
\def\epsilon{\varepsilon}

\newcommand{\md[1]}{{\rm Mod}({#1})}

\newcommand{\To}{{\makebox[2em]{\rightarrowfill}}}
\def\Op{{\rm Op}}
\def\Xsa{X_{sa}}
\def\dist{{\rm dist}}

\newcommand{\Indo[1]}{{\rm Ind}({#1})}
\newcommand{\II[1]}{{\rm I}({#1})}

\def\BBP{{\mathbb P}}

\newcommand{\ba}{\begin{array}}
\newcommand{\ea}{\end{array}}

\def\shb{\mathcal{B}}
\def\shc{\mathcal{C}}
\def\shd{\mathcal{D}}

\def\shi{\mathcal{I}}

\def\shm{\mathcal{M}}

\def\sho{\mathcal{O}}

\def\shs{\mathcal{S}}
\def\sht{\mathcal{T}}
\def\shu{\mathcal{U}}

\newcommand{\C}{\mathbb{C}}

\newcommand{\R}{\mathbb{R}}
\newcommand{\Q}{\mathbb{Q}}
\newcommand{\Z}{\mathbb{Z}}

\renewcommand{\to}[1][]{\xrightarrow[#1]{}}
\newcommand{\from}[1][]{\xleftarrow[#1]{}}
\newcommand{\isoto}[1][]{\xrightarrow[#1]{\sim}}

\DeclareMathOperator{\ori}{or}
\newcommand{\mdco[1]}{{\rm Mod_{coh}}({#1})}
\newcommand{\mdc[1]}{{\rm Mod}^c({#1})}
\newcommand{\ssubset}{\subset\subset}
\newcommand{\tens}{\otimes}

\newcommand{\rc[1]}{{\Brc}({#1})}
\newcommand{\rcc[1]}{{\Brc}^c({#1})}

\newcommand{\eim}[1]{{#1}_{!}}
\newcommand{\oim}[1]{{#1}_*}
\newcommand{\opb}[1]{#1^{-1}}
\newcommand{\eeim}[1]{{#1}_{!!}}
\renewcommand{\hom}[1][]{{\mathcal{H}om}_{\raise1.5ex\hbox to.1em{}#1}}
\newcommand{\ihom}[1][]{{\shi hom}_{\raise1.5ex\hbox to.1em{}#1}}
\newcommand{\reeim}[1]{R{#1}_{!!}}
\newcommand{\roim}[1]{{R#1}_*}
\newcommand{\reim}[1]{{R#1}_!}
\newcommand{\epb}[1]{#1^{!}}
\newcommand{\rhom}[1][]{{R\mathcal{H}om}_{\raise1.5ex\hbox to.1em{}#1}}
\newcommand{\thom}[1][]{{T\mathcal{H}om}_{\raise1.5ex\hbox to.1em{}#1}}
\newcommand{\tmuhom}[1][]{{T\mu\mathcal{H}om}_{\raise1.5ex\hbox to.1em{}#1}}
\newcommand{\tPhi}{\tilde{\Phi}}

\newcommand{\sect}{\Gamma}
\newcommand{\rsect}{R\Gamma}
\newcommand{\risect}{R{\rm I}\Gamma}
\newcommand{\isect}{{\rm I}\Gamma}

\newcommand{\Hom}[1][]{\mathrm{Hom}_{\raise1.5ex\hbox to.1em{}#1}}
\newcommand{\RHom}[1][]{\mathrm{RHom}_{\raise1.5ex\hbox to.1em{}#1}}
\def\wtens{\stackrel{{\rm w}}{\otimes}} 

\newcommand{\IIrc[1]}{{\rm I_{\R-c}}({#1})}
\newcommand{\IRc[1]}{{\rm I}\R{\rm-c}({#1})}
\newcommand{\IRC}{{\rm I}\,\R{\rm-c}}
\newcommand{\wRC}{{\rm w-}\R{\rm-c}}
\newcommand{\wCC}{{\rm w-}\C{\rm-c}}
\newcommand{\RC}{\R{\rm-c}}
\newcommand{\CC}{\C{\rm-c}}

\def\car{{\rm char}}

\theoremstyle{plain}

\newtheorem{theorem}{Theorem}[section]
\newtheorem{proposition}[theorem]{Proposition}
\newtheorem{lemma}[theorem]{Lemma}

\newtheorem{conjecture}[theorem]{Conjecture}

\theoremstyle{definition}

\newtheorem{definition}[theorem]{Definition}

\newtheorem{example}[theorem]{Example}

\newtheorem{remark}[theorem]{Remark}

\newenvironment{nnum}{
  \begin{enumerate}
  \itemsep=0pt
  
  }
  {\end{enumerate}}

\newenvironment{anum}{
  \begin{enumerate}
  \itemsep=0pt
  
  }
  {\end{enumerate}}

\newcommand{\bnum}{\begin{nnum}}
\newcommand{\enum}{\end{nnum}}
\newcommand{\banum}{\begin{anum}}
\newcommand{\eanum}{\end{anum}}
\newcommand{\ol}{\overline}
\def\dddt{{\raise-.3em\hbox{$\big\cdot$}}}

\newcommand{\cl}{\colon}
\newcommand{\hs}[1]{\hspace*{#1}}
\newcommand{\Real}{\operatorname{Re}}

\usepackage[all]{xy}
\CompileMatrices
\numberwithin{equation}{section}

\author{Masaki Kashiwara and Pierre Schapira}

\title{Microlocal study of ind-sheaves I:\\
micro-support and regularity}
\date{June 2001}

\begin{document}

\maketitle

\begin{abstract}
We define the notions of micro-support and regularity
for ind-sheaves, and prove their invariance
by contact transformations.
We apply the results
to the ind-sheaves of temperate holomorphic
solutions of $\shd$-modules.
We  prove that the micro-support of such
an ind-sheaf is the 
characteristic variety of the corresponding $\shd$-module 
and that the ind-sheaf is regular if the $\shd$-module is regular holonomic.
\end{abstract}
\renewcommand{\thefootnote}{}

\tableofcontents

\section{Introduction}
Recall that a system of linear partial
differential equations on a complex
manifold $X$ is the data of a coherent module $\shm$ over the sheaf of
rings $\shd_X$ of holomorphic differential operators. 
Let $F$ be a complex of sheaves on $X$ with $\R$-constructible cohomologies
(one says an $\R$-constructible sheaf, for short).
The complex of ``generalized functions'' associated with $F$ is
described by the complex $\rhom(F,\sho_X)$, 
and the complex of solutions of $\shm$
with values in this complex is described by the complex
$$
\rhom[\shd_X](\shm,\rhom(F,\sho_X)).
$$

One may also microlocalize the problem by replacing
$\rhom(F,\sho_X)$ with $\mu hom(F,\sho_X)$. In \cite{K-S1} one shows that
most of the properties of this complex, especially those related to
propagation or Cauchy problem,  are encoded in two geometric objects,
both living in the cotangent bundle $T^*X$, the characteristic variety
of the system
$\shm$, denoted by $\car(\shm)$, and the micro-support of $F$, denoted
by $SS(F)$.

The complex $\rhom(F,\sho_X)$ allows us to treat various situations.
For example if $M$ is a real analytic 
manifold and $X$ is a complexification of $M$,
by taking as $F$ the dual $D'(\C_M)$ of the constant sheaf on $M$, 
one obtains the sheaf $\shb_M$ of Sato's hyperfunctions. If $Z$ is a complex
analytic hypersurface of $X$ and $F=\C_Z[-1]$ is the (shifted)
constant sheaf on $Z$, one obtains 
the sheaf of holomorphic functions with singularities on $Z$.
However, the complex $\rhom(F,\sho_X)$ does not allow us to 
treat sheaves associated with holomorphic functions  with growth 
conditions. So far this difficulty was overcome in two cases, 
the temperate case including Schwartz's distributions and meromorphic
functions with poles on $Z$ and the dual case
including $C^\infty$-functions and the formal completion of $\sho_X$ along
$Z$. The method was two construct specific functors, the functor 
$\thom$ of \cite{K1} and the functor $\wtens$ of \cite{K-S2}.

There is a more radical method, which consists in replacing the
too narrow framework of sheaves by that of ind-sheaves, as explained
in \cite{K-S3}. For example,
the presheaf of holomorphic
temperate functions on a complex manifold $X$ (which, to a subanalytic 
open subset of $X$,
associates the space of holomorphic functions with temperate growth
at the boundary) is clearly not a sheaf.
However it makes sense as an
object (denoted by ${\sho}^t_X$) of the derived category of ind-sheaves on $X$.
Then it is natural to ask if the microlocal theory of sheaves, in
particular the theory of micro-support, applies in this general
setting. 

In this paper we give the definition and the elementary properties of the
micro-support of ind-sheaves as well as the notion of regularity. 

We prove in particular that the micro-support $SS(\cdot)$ and 
the regular micro-support
$SS_{reg}(\cdot)$ of ind-sheaves behave naturally with respect to
distinguished triangles and 
that these micro-supports are invariant 
by ``quantized contact transformations'' 
(in the framework of sheaf theory, as explained in \cite{K-S1}).

When $X$ is a complex manifold and $\shm$ is a coherent
$\shd_X$-module, we study the ind-sheaf   
$Sol^t(\shm):=\rhom[\shd_X](\shm,{\sho}^t_X)$. We prove that
\bnum
\item
$SS(Sol^t(\shm))=\car(\shm)$,
\item  if $\shm$ is holonomic, $Sol^t(\shm)$ is regular if 
$\shm$ is regular holonomic.
\enum
Finally, we treat an
example: we calculate  the
ind-sheaf of the temperate holomorphic solutions 
of an irregular differential equation.

This paper is the first one of a series. In Part II, we shall
introduce the microlocalization functor for ind-sheaves, and in Part
III we shall study the functorial behavior of micro-supports.

\section{Notations and review}\label{section:notation}
We will mainly follow the notations in \cite{K-S1} and \cite{K-S3}.

\vspace{0.5 cm}
\noindent
{\bf Geometry.}

\noindent 
In this paper, all manifolds will be real analytic 
(sometimes, complex analytic).
Let $X$ be a manifold. One denotes by 
$\tau\colon TX\to X$ the tangent bundle to $X$ and
by $\pi\colon T^*X\to X$ the
 cotangent bundle. One denotes by $a\cl T^*X\to T^*X$ the  antipodal
map. If $S\subset T^*X$, 
one denotes by $\dot{S}$ the set $S\setminus T^*_XX$, 
and one denotes by $S^a$ the image of $S$ by the antipodal
map. In particular, $\dot{T}^*X=T^*X\setminus X$, the set $T^*X$ with the
zero-section removed. One denotes by 
$\dot{\pi}\colon\dot{T}^*X\to X$  the projection.

For a smooth submanifold $Y$ of $X$, $T_YX$ denotes
 the normal bundle to $Y$ and $T^*_YX$ the conormal bundle. 
In particular, $T^*_X X$ is identified with $X$, the zero-section. 

For a submanifold $Y$ of $X$ and a subset  $S$ of $X$, 
we denote by $C_Y(S)$ the Whitney normal cone to $S$
along $Y$,  a conic subset of $T_YX$. 

If $S$ is  a locally closed subset of
 $T^*X$, we say that $S$ is $\R^+$-conic (or simply ``conic'', for
 short) if it is locally invariant under 
the action of $\R^+$. If $S$ is smooth, this is equivalent to
saying that the Euler vector field on $T^*X$ is tangent to $S$.

Let $f\colon X\to Y$ be a morphism of real manifolds. 
One has two natural maps

\eq\label{eq:fdandfpi}
T^*X\from[f_d]X\times_Y T^*Y\to[f_\pi] T^*Y
\eneq
(In \cite{K-S1}, $f_d$ is denoted by ${}^{t}f'$.)
We denote by $q_1$ and $q_2$ 
the first and second projections defined
on $X\times Y$. 

\vspace{0.5 cm}
\noindent
{\bf Sheaves.} 

\noindent
Let $k$ be a field. We denote by $\md[k_X]$ the abelian category 
of sheaves of $k$-vector spaces and by $D^b(k_X)$ its
 bounded derived category.

We denote by $\rc[k_X]$ the abelian category of $\R$-constructible
sheaves of $k$-vector spaces on $X$, and by 
$D^b_{\RC}(k_X)$ (resp. $D^b_{\wRC}(k_X)$) 
the full triangulated
subcategory of $D^b(k_X)$ consisting of objects with $\R$-constructible 
(resp. weakly $\R$-constructible) cohomology.
On a complex manifold,
one defines similarly the categories $D^b_{\CC}(k_X)$
and $D^b_{\wCC}(k_X)$ of $\C$-constructible and weakly $\C$-constructible
sheaves.

If $Z$ is  a locally closed subset of $X$ and if $F$ is  a sheaf on
$X$, recall that $F_Z$ is a sheaf on $X$
such that $F_Z\vert_Z\simeq F\vert_Z$ and
$F_Z\vert_{X\setminus Z}\simeq 0$. One writes $k_{XZ}$ instead of
$(k_X)_Z$ and one sometimes writes $k_Z$ instead of
$k_{XZ}$.

If $f\cl X\to Y$ is a morphism of manifolds, one denotes by
$\omega_{X/Y}$ the relative dualizing complex on $X$ and if
$Y=\{\rmpt\}$ one simply denotes it by $\omega_X$.
Recall that
$$\omega_X\simeq {\ori}_X [\dim_\R X]$$
where $\ori_X$ is the orientation sheaf and $\dim_R X$ is the dimension of
$X$ as a real manifold. 
We denote by $D'_X$ and $D_X$ the duality functors on $D^b(k_X)$, defined by 
$$ D'_X(F)=\rhom(F,k_X),\,\,  D_X(F)=\rhom(F,\omega_X).$$

If $F$ is an object of $D^b(k_X)$, $SS(F)$ denotes its micro-support, a closed
conic involutive subset of $T^*X$. For an open subset $U$ of $T^*X$,
one denotes by $D^b(k_X;U)$ the localization of the category
$D^b(k_X)$ with respect to the triangulated 
subcategory consisting of sheaves $F$
such that $SS(F)\cap U=\emptyset$.

We shall also use the functor $\mu hom$ as well as the operation 
$\widehat{+}$ and refer to loc. cit.
for details.

\vspace{0.5 cm}
\noindent
{\bf $\sho$ and $\shd$}

\noindent
On a complex manifold $X$ we consider the structural sheaf  $\sho_X$  of
holomorphic functions and the sheaf $\shd_X$ 
of linear holomorphic differential operators of finite order.

We denote by 
$\mdco[\shd_X]$  
the abelian category of  coherent $\shd_X$-modules.
We denote by $D^b(\shd_X)$ the bounded derived category of left
$\shd_X$-modules and by $D^b_{\rm coh}(\shd_X)$ 
(resp. $D^b_{\rm hol}(\shd_X))$, $D^b_{\rm rh}(\shd_X))$ 
its full triangulated
category consisting of objects with coherent cohomologies
(resp. holonomic cohomologies, regular holonomic cohomologies). 


\vspace{0.5 cm}
\noindent
{\bf Categories.} 
In this paper, we shall work in a given universe $\shu$, and a
category means a $\shu$-category. If $\shc$ is a category,
$\shcw$ denotes the category of functors from $\shc^{\rop}$ to
$\Set$. The category $\shcw$ admits inductive limits, however, in case
$\shc$ also admits inductive limits, the Yoneda functor 
$h^\wedge\cl \shc\to\shcw$ does not 
commute with such limits. Hence, one denotes by  $\indlim$ the
inductive limit in $\shc$ and by $\inddlim$ the inductive limit in $\shcw$.

One denotes by
$\Indo[\shc]$ the category of ind-objects of $\shc$, that is the full
subcategory of $\shcw$ consisting of objects $F$ such that there
exist a small filtrant category $I$ and a functor $\alpha\cl I\to \shc$,
with 
$$
F\simeq\inddlim\alpha, \mbox{ i.e., }
F\simeq\inddlim[i\in I]F_i, \mbox{ with }F_i\in\shc.
$$
The category $\shc$ is considered as a full subcategory of $\Indo[\shc]$.

If $\phi\cl \shc\to\shc'$ is  a functor, it defines a functor 
$I\phi\cl \Indo[\shc]\to\Indo[\shc']$ which commutes with
$\inddlim$. 

If $\shc$ is an additive category, we denote by $C(\shc)$ the category
of complexes in $\shc$ and by $K(\shc)$ the associated homotopy
category. If $\shc$ is abelian, one denotes by $D(\shc)$ its derived
category. One defines as usual the full subcategories 
$C^*(\shc), K^*(\shc), D^*(\shc)$, with $*=+,-,b$. One denotes by 
$Q$ the localization functor:
\eqn
Q\cl K^*(\shc)\to D^*(\shc).
\eneqn 
We keep the same notation $Q$ to
denote the composition $C^*(\shc)\to K^*(\shc)\to D^*(\shc)$.

One denotes by $C^{[a,b]}(\shc)$ the full subcategory 
of  $C(\shc)$ consisting of objects
$F^\bullet$ satisfying $F^i=0$ for $i\notin [a,b]$.
If $a,b\in\Z$ with $a\leq b$, there is a natural isomorphism
$$\Indo[C^{[a,b]}(\shc)]\isoto C^{[a,b]}(\indc).$$

\vspace{0.5 cm}
\noindent
{\bf Ind-sheaves.}
Here, $X$ is a Hausdorff locally compact space with a countable
base of open sets
and $k$ is a field.
One denotes by $\II[k_X]$ the abelian category of ind-sheaves of
$k$-vector spaces on $X$, that is, $\II[k_X]= \Indo[{\mdc[k_X]}]$, the
category of ind-objects of the category $\mdc[k_X]$ of sheaves with
compact support on $X$. 
We denote by $D^b(\II[k_X])$ the bounded derived category of
$\II[k_X]$.

There is a natural fully faithful exact functor 
\eqn
&&\iota_X\cl \md[k_X]\to \II[k_X],\\
&&F\mapsto \inddlim[U\ssubset X]F_U \mbox{ ($U$ open)}.
\eneqn
Most of the time, we shall not write this functor and identify
$\md[k_X]$ with a full abelian subcategory of $\II[k_X]$ and $D^b(k_X)$
with a full triangulated subcategory of $D^b(\II[k_X])$.

The category $\II[k_X]$ admits
an internal hom denoted by $\ihom$ and this functor admits a left adjoint,
denoted by $\tens$. If $F\simeq\inddlim[i]F_i$ and $G\simeq\inddlim[j]G_j$, then
\eqn
\ihom(G,F)&\simeq&\prolim[j]\inddlim[i]\hom(G_j,F_i)\\
G\tens F &\simeq&\inddlim[i]\inddlim[j](G_j\tens F_i).
\eneqn

The functor $\iota_X$ admits a left adjoint
\eqn
&&\alpha_X\cl \II[k_X]\to\md[k_X],
\eneqn 
To 
$F=\inddlim[i\in I]F_i$,  this functor 
associates $\alpha_X(F)=\indlim[i\in I]F_i$.
This functor also admits a left adjoint 
\eqn
&&\beta_X\cl \md[k_X]\to \II[k_X], 
\eneqn
and both functors $\alpha_X$ and $\beta_X$ are exact. 
The functor $\beta_X$ is not so easy to describe. For example, 
for an open subset $U$ and a closed subset $Z$, one has;
\eqn
&&\beta_X(k_{XU})\simeq\inddlim[V\ssubset U]k_{XV}\mbox{ ($V$ open),}\\
&&\beta_X(k_{XZ})\simeq\inddlim[Z\subset V]k_{X\ol{V}}\mbox{ ($V$ open).}
\eneqn
One sets
\eqn
\hom(G,F)=\alpha_X\ihom(G,F)\in \md[k_X].
\eneqn
One has $$\Hom[{\II[k_X]}](G,F)=\sect(X;\hom(G,F)).$$
The functors $\ihom$ and $\hom$ are left exact and admit right derived
functors $R\ihom$ and $\rhom$.

Let $f\cl X\to Y$ be a morphism of topological spaces ($Y$
satisfies the same assumptions as $X$). 
There are natural functors
\eqn
&&\opb{f}\cl  \II[k_Y])\to\II[k_X]\\
&&\oim{f}\cl  \II[k_X]\to \II[k_Y]\\
&&\eeim{f}\cl  \II[k_X]\to \II[k_Y].
\eneqn
The proper direct image functor is denoted by $\eeim{f}$ instead of
$\eim{f}$ because it does not commute with $\iota$, that is 
$\iota_Y\eim{f}\neq\eeim{f}\iota_X$ in general..

These functors induce derived functors, and moreover the functor 
$\reeim{f}$ admits a right adjoint denoted by $\epb{f}$:
\eqn
&&\opb{f}\cl  D^b\II[k_Y])\to D^b(\II[k_X]),\\
&&\roim{f}\cl  D^b(\II[k_X])\to D^b(\II[k_Y]),\\
&&\reeim{f}\cl  D^b(\II[k_X])\to D^b(\II[k_Y]),\\
&&\epb{f}\cl  D^b(\II[k_Y])\to D^b(\II[k_X]).
\eneqn
Let $a_X\cl X\to\{\rmpt\}$ denote the canonical map. We also introduce a
notation. We set
\eqn
\isect(X;\cdot)&=&\oim{a_X}(\cdot),\\
\risect(X;\cdot)&=&\roim{a_X}(\cdot).
\eneqn

\vspace{0.5 cm}
\noindent
{\bf Ind-sheaves on real manifolds.} 
Let $X$ be a real analytic manifold.
Among all ind-sheaves, there are those which are ind-objects of 
the category of $\R$-constructible sheaves, and  we shall
encounter them in our applications. 

We denote by $\rcc[k_X]$ the full
abelian subcategory of $\rc[k_X]$ consisting of 
$\R$-constructible sheaves with compact support.
We set 
$$\IRc[k_X]=\Indo[{\rcc[k_X]}]$$
and denote by $D^b_{\IRC}(\II[k_X])$ the full 
subcategory of $D^b(\II[k_X])$ consisting of objects with cohomology
in $\IRc[k_X]$. 
(Note that in \cite{K-S3}, $\IRc[k_X]$ was denoted by $\IIrc[k_X]$.)

\begin{theorem}\label{eq:dbirceqv}
The natural functor 
$D^b(\IRc[k_X])\to D^b_{\IRC}(\II[k_X])$
is an equivalence. 
\end{theorem}
There is an alternative construction of $\IRc[k_X]$, using
Grothendieck topologies. Denote by $\Op_X$ the
category of open subsets of $X$ (the morphisms $U\to V$ are the inclusions),
and by by $\Op_{\Xsa}$ its full subcategory consisting 
of open subanalytic subsets of $X$.
One endows this category with a Grothendieck
topology by deciding that a family $\{U_i\}_i$ in $\Op_{\Xsa}$ 
is a covering of $U\in\Op_{\Xsa}$ if for any compact subset $K$ of $X$,
there exists a finite subfamily which covers $U\cap K$. In other words,
we consider families which are locally finite in $X$. One
denotes by $\Xsa$ the site defined by this topology.

Sheaves on $\Xsa$ are easy to construct. Indeed, 
consider a
presheaf $F$ of $k$-vector spaces defined on the subcategory
$\Op^c_{\Xsa}$ of relatively
compact open subanalytic subsets of $X$ and assume that the sequence
$$
0\to F(U\cup V)\to F(U)\oplus F(V)\to F(U\cap V)$$
is exact for any $U$ and $V$ in $\Op^c_{\Xsa}$. 
Then there exists a unique sheaf $\tilde F$ on $\Xsa$ 
such that $\tilde F(U)\simeq F(U)$ for all $U\in\Op^c_{\Xsa}$.
Sheaves on $\Xsa$ define naturally ind-sheaves on $X$. Indeed:

\begin{theorem}\label{th:suba}
There is a natural equivalence of abelian categories
$$\IRc[k_X]\isoto\md[k_{\Xsa}],$$
given by
$$\IRc[k_X]\ni F\mapsto
\bigl(\Op^c_{\Xsa}\ni U\mapsto\Hom[{\IRc[k_X]}](k_U,F)\bigr).$$
\end{theorem}
As usual, we denote by 
$\shc_X^\infty$  the sheaf of complex-valued functions
of class $\shc^\infty$, by $\shd b_X$ 
(resp.\ $\shb_X$) 
the sheaf of Schwartz's distributions (resp.\ Sato's hyperfunctions),
and by $\shd_X$ the sheaf of analytic finite-order differential operators.

Let $U$ be an open subset of $X$. One sets
$\shc^{\infty}_X(U)=\sect(U;\shc^{\infty}_X)$. 

\begin{definition}
Let $f\in \shc^{\infty}_X(U)$. One says that $f$ has
{\it  polynomial growth} at $p\in X$
if it satisfies the following condition.
For a local coordinate system
$(x_1,\dots,x_n)$ around $p$, there exist
a sufficiently small compact neighborhood $K$ of $p$
and a positive integer $N$
such that
\begin{eqnarray}
&\sup_{x\in K\cap U}\big(\dist(x,K\setminus U)\big)^N\vert f(x)\vert
<\infty\,.&
\end{eqnarray}

It is obvious that $f$ has polynomial growth at any point of $U$.
We say that $f$ is {\it tempered} at $p$ if all its derivatives
have polynomial growth at $p$. We say that $f$ is tempered 
if it is tempered at any point.
\end{definition}

For an open subanalytic set $U$ in $X$,
denote by $\shc_X^{\infty,t}(U)$ the subspace
of $\shc^{\infty}_X(U)$ consisting of tempered functions. 
Denote by $\shd b_X^t(U)$ the space of tempered distributions on $U$, 
defined by the exact sequence
$$0\to\sect_{X\setminus U}(X;\shd b_X)\to\sect(X;\shd b_X)
\to\shd b_X^t(U)\to 0.$$
It follows from the results of Lojasiewicz \cite{Lo} that 
$U\mapsto \shc^{\infty}_X(U)$ and $U\mapsto \shd b_X^t(U)$
are sheaves on the subanalytic site $\Xsa$, hence define ind-sheaves.
\begin{definition}
We call $\shc^{\infty,t}_X$ 
(resp. $\shd b_X^t$)
the ind-sheaf of tempered $\shc^{\infty}$-functions 
(resp. tempered distributions).
\end{definition}
One can also define the ind-sheaf of Whitney
$\shc^{\infty}$-functions, but we shall not recall here its
construction.
These ind-sheaves are well-defined in the category
$\md[\beta_X\shd_X]$. Roughly speaking, it means that if 
$P$ is a differential operator defined on the closure $\bar U$ of an
open subset $U$, then it acts on $\shc^{\infty,t}_X(U)$ and $\shd b_X^t(U)$.

Let now $X$ be  a complex manifold.
We denote by 
$\overline{X}$ the complex conjugate manifold and by $X^\R$ the underlying 
real analytic manifold, identified with the diagonal of $X\times \overline{X}$.
We denote by $\shd_X$ the sheaf
of rings of finite-order holomorphic differential operators, not
to be confused with $\shd_{X^\R}$. We set
\eqn
\sho_X^t&:=&
  R\ihom[\beta\shd_{\overline{X}}](\beta\sho_{\overline{X}},\shd b^t_{X^\R})
\eneqn
One can prove that the natural morphism
\eqn
R\ihom[\beta\shd_{\overline{X}}]
(\beta\sho_{\overline{X}},\shc^{\infty,t}_{X^\R})\to
R\ihom[\beta\shd_{\overline{X}}](\beta\sho_{\overline{X}},\shd b^t_{X^\R})
\eneqn
is an isomorphism.
One calls $\sho_X^t$ the ind-sheaf of tempered holomorphic
functions. One shall be aware that in fact, $\sho_X^t$ is not an
ind-sheaf but an object of the derived category $D^b(\II[\C_X])$, or
better, of $D^b(\beta_X\shd_X)$. It
is not concentrated in degree $0$ as soon as $\dim X>1$.

Let $G\in D^b_{\R-c}(\C_X)$.
It follows from the construction of $\sho_X^t$ that:
\eqn
R\hom(G,\sho_X^t)\simeq \thom(G,\sho_X),
\eneqn
where $\thom(\cdot,\sho_X)$ denotes the functor of temperate
cohomology of \cite{K1} (see also \cite{K-S2} for a detailed
construction and \cite{An} for its microlocalization).

\section{Complements of homological algebra}\label{section:homal}

The results of this section are extracted from \cite{K-S4}.
Let $\shc$ denote a small abelian category.
We shall study some links between the derived category
$D^b(\indc)$ and the category $\Indo[D^b(\shc)]$. 

We define the functor 
$J\cl D^b(\indc)\to (D^b(\shc))^\wedge$
by setting for $F\in D^b(\indc)$ and $G\in D^b(\shc)$

\begin{equation}\label{eq:J(X)1}
J(F)(G)=\Hom[D^b(\indc)](G,F).
\end{equation}
\begin{theorem}\label{th:indcab}
\bnum
\item
The functor $J$ takes its values in $\Indo[D^b(\shc)]$.
\item
Consider a small and filtrant category $I$, integers $a\le b$
and a functor $I\to C^{[a,b]}(\shc)$, $i\mapsto F_i$. 
If $F\in D^b(\indc)$, $F\simeq Q(\smash{\inddlim[i]F_i})$ and 
$G\in D^b(\shc)$, then:
\banum
\item
$J(F)\simeq \inddlim[i]Q(F_i)$,
\item
$\Hom[{D^b(\indc)}](G,F)\simeq\indlim[i]\Hom[D^b(\shc)](G,F_i).$
\eanum
\item
For each $k\in\Z$, the diagram below commutes.
$$
\xymatrix{
D^b({\indc})\ar[rr]_-J\ar[rd]_-{H^k}& & {\Indo[D^b(\shc)]}
 \ar[ld]^-{IH^k}\\
  & {\indc} & 
}$$
\enum
\end{theorem}

\begin{lemma}\label{le:isoinind}
Assume that $\shc$  has  finite homological dimension.
Let $\phi\cl X\to Y$ be a morphism in $\Indo[D^b(\shc)]$
and assume that $\phi$ induces an isomorphism 
$IH^k(\phi)\cl IH^k(X)\isoto IH^k(Y)$ for every $k\in \Z$. 
Then $\phi$ is an isomorphism.
\end{lemma}

\begin{theorem}\label{th:indcab2}
Let $\psi\cl D^b(\Indo[\shc])\to D^b(\Indo[\shc'])$ be a triangulated functor
which satisfies: if $F\in D^b(\Indo[\shc])$, $F\simeq Q(\inddlim[i]F_i)$ with 
$F_i\in C^{[a,b]}(\shc)$, 
then $H^k\psi(F)\simeq \inddlim[i]H^k\psi(Q(F_i))$. 
Assume moreover that the homological dimension of $\shc'$ is finite.
Then there exists
a unique functor $J\psi\cl \Indo[D^b(\shc)]\to \Indo[D^b(\shc')]$ which
commutes with $\inddlim$  and
such that the diagram below commutes:
$$
\xymatrix{
{D^b(\indc)}\ar[r]_-{\psi}\ar[d]_-{J}&D^b(\Indo[\shc'])\ar[d]_-{J}\\
{\Indo[D^b(\shc)]}\ar[r]^-{J\psi}    &{\ \Indo[D^b(\shc')]\,.}
}$$
\end{theorem}

\begin{remark}
The functor $J\cl D^b(\indc)\to\Indo[D^b(\shc)]$ is neither full nor
faithful. 
Indeed, let $\shc=\mdc[k_X]$ and let $F\in \md[k_X]$ considered as a
full subcategory of $\II[k_X]$. Then 
$$\Hom[D^b({\II[k_X]})](k_X,F[n])\simeq H^n(X;F).$$
On the other hand,
$$\Hom[{\Indo[{D^b(\mdc[k_X])}]}](J(k_X),J(F[n]))\simeq\prolim[U\ssubset X]H^n(U;F).$$
\end{remark}

Let $\sht$ be a full triangulated subcategory $D^b(\shc)$. One
identifies $\Indo[\sht]$ with a full subcategory of $\Indo[D^b(\shc)]$.

Let $F\in D^b(\indc)$. Let us denote by $\sht_F$ the category
of arrows $G\to F$ in $D^b(\indc)$ with $G\in\sht$.
The category  $\sht_F$ is filtrant.

\begin{lemma}\label{le:Jshdcohom}
For $F\in D^b(\indc)$,
the conditions below are equivalent.
\bnum
\item
$J(F)\in\Indo[\sht]$,
\item
for each $k\in\Z$, one has $H^k(F)\simeq \inddlim[G\to F\in\sht_F]H^k(G)$.
\enum
\end{lemma}

\begin{definition}
Let $\sht$ be a full triangulated subcategory of $D^b(\shc)$. One
denotes by $\opb{J}\Indo[\sht]$ the full subcategory 
of $D^b(\indc)$ consisting of
objects $F\in D^b(\indc)$ such that $J(F)\in\Indo[\sht]$.
\end{definition}

\begin{proposition}\label{th:ishttriang}
The category $\opb{J}\Indo[\sht]$ is a 
triangulated subcategory of $D^b(\indc)$.
\end{proposition}

\bigskip
We will apply these results to the category $\II[k_X]=\Ind(\mdc[k_X])$.
Hence $J$ is the functor:
$$J\cl  D^b(\II[k_X])\To \Indo[{D^b(\mdc[k_X])}].$$
By the definition one has
$$\mbox{$J(F)\simeq\inddlim[{U\ssubset X}] J(F_U)$ 
\quad for any $F\in D^b(\II[k_X])$.}$$

%
%
As a corollary of Theorem \ref{th:indcab2}, one gets:
\begin{proposition}\label{th:JtensandJihom} 
For $G\in D^b(k_X)$ and $F\in D^b(\II[k_X])$,
assume that $J(F)\simeq\smash{\inddlim[i]J(F_i)}$ with $F_i\in D^b(k_X)$.
Then there are natural isomorphisms:
\eq
J(G\tens F)&\simeq&\inddlim[i]J(G\tens F_i),\label{eq:jtens}\\
J(R\ihom(G,F))&\simeq &\inddlim[i]J(R\ihom(G,F_i)).\label{eq:jihom}
\eneq
\end{proposition}

\section{Micro-support and regularity}\label{section:micro-support}

Let $\gamma$ be a closed convex proper
cone in an affine space $X$. One denotes by $\gamma^{\circ}$ its polar cone,
\eqn
&&\gamma^{\circ}=\{\xi\in X^*;\langle x,\xi\rangle\geq 0
\mbox{ for all $x\in\gamma\}$.}
\eneqn
Let $W\subset X$ be an open subset. We introduce the
functor $\Phi_{\gamma,W}\cl  D^b(\II[k_X])\to D^b(\II[k_X])$ as follows.
Denote by $q_1,q_2\cl X\times X\to X$ the first and second projections
and  denote by $s\cl X\times X\to X$
the map $(x,y)\mapsto x-y$. One sets
\eqn
&&\Phi_{\gamma,W}(F)=\reeim{q_1}(k_{\opb{s}\gamma\cap\opb{q_1}W\cap\opb{q_2}W}
\tens\opb{q_2}F).
\eneqn
One writes $\Phi_{\gamma}$ instead of $\Phi_{\gamma,X}$.
Define the functor ${\Phi}^-_{\gamma,W}$ by replacing 
the kernel $k_{\opb{s}\gamma\cap\opb{q_1}W\cap\opb{q_2}W}$ with 
the complex $k_{\opb{s}\gamma\cap\opb{q_1}W\cap\opb{q_2}W}\to k_{s^{-1}(0)}$
in which $k_{s^{-1}(0)}$ is situated in degree $0$.
We have a distinguished triangle in $D^b(\II[k_X])$
\eqn
\Phi_{\gamma,W}(F)\to F\to {\Phi}^-_{\gamma,W}(F)\xrightarrow{+1}.
\eneqn

Note that if $F\in D^b(k_X)$, then 
\eqn
\left\{ \begin{array}{l}
\supp(\Phi_{\gamma,W}(F))\subset \ol{W},\\
\Phi_\gamma(F)\to F \mbox{ is an isomorphism on }X\times\Int{\gamma}^\circ,\\
SS(\Phi_\gamma(F))\subset X\times{\gamma}^\circ.\\
SS({\Phi}^-_{\gamma,W}(F))\bigcap W\times\Int{\gamma}^\circ=\emptyset
\end{array}\right.
\eneqn

\begin{lemma}\label{le:ssf}
Let $F\in D^b(\II[k_X])$ and let $p\in T^*X$. 
The conditions {\rm (1a)--(4b)} below are all equivalent.
Moreover, if $F\in D^b_{\IRC}(\II[k_X])$, these conditions are
equivalent to {\rm (5a)}.

\bnum
\item[{\rm (1a)}]
Assume that for a small and
filtrant category $I$, integers $a\leq b$
and a functor $I\to C^{[a,b]}(\md[k_X])$, $i\mapsto F_i$ one has
$F\simeq Q(\inddlim[i\in I] F_i)$.
Then there exists a conic open neighborhood $U$ of $p$ in $T^*X$ such that
for any $i\in I$ there exists a morphism $i\to j$ in $I$ which induces
the zero-morphism $0:F_i\to F_j$ in $D^b(k_X;U)$.
\item[{\rm (1b)}]
There exist a conic open neighborhood $U$ of $p$ in $T^*X$, 
a small and
filtrant category $I$, integers $a\le b$ 
and a functor $I\to C^{[a,b]}(\md[k_X])$, $i\mapsto F_i$, such that 
$SS(F_i)\cap U=\emptyset$ and $F\simeq Q(\smash{\inddlim[i] F_i})$ in a
neighborhood of $\pi(p)$.
\item[{\rm (2a)}]
Assume that for a small and
filtrant category $I$, integers $a\leq b$
and a functor $I\to D^{[a,b]}(k_X)$, $i\mapsto F_i$ one has
$J(F)\simeq \inddlim[i\in I] J(F_i)$.
Then there exists a conic open neighborhood $U$ 
of $p$ in $T^*X$ such that
for any $i\in I$ there exists a morphism $i\to j$ in $I$ which induces
the zero-morphism $0:F_i\to F_j$ in $D^b(k_X;U)$.
\item[{\rm (2b)}]
There exist a conic open neighborhood $U$ of $p$ in $T^*X$, 
a small and filtrant category $I$, integers $a\leq b$,
a functor $I\to D^b(k_X)$, $i\mapsto F_i$ and $F'$
isomorphic to $F$ in neighborhood of $\pi(p)$ such that 
$SS(F_i)\cap U=\emptyset$ and $J(F')\simeq\inddlim[i]J(F_i)$.
\item[{\rm (3a)}] 
There exists a conic open neighborhood $U$ of $p$ in $T^*X$ such that for
any $G\in D^b(k_X)$ with $\supp (G)\ssubset \pi(U)$, 
$SS(G)\subset U\cup T^*_XX$, one has
$\Hom[{D^b(\II[k_X])}](G,F)$ $=0$.
\item[{\rm (3b)}] 
There exists a conic open neighborhood $U$ of $p$ in $T^*X$ such that for
any $G\in D^b(k_X)$ with $\supp (G)\ssubset \pi(U)$, 
$SS(G)\subset {U}^a\cup T^*_XX$, one has 
$\risect(X;G\tens F)=0$.
\enum
Assume now that $X$ is an affine space
and let $p=(x_0;\xi_0)$.
\bnum
\item[{\rm (4a)}]
There exist a relatively compact open neighborhood $W$ of $x_0$ and 
a closed convex proper cone $\gamma$ with $\xi_0\in\Int \gamma^\circ$
such that $\Phi_{\gamma,W}(F)\simeq 0$.
\item[{\rm (4b)}] There exist $F'\in  D^b(\II[k_X])$ with $F'\simeq F$ in a
neighborhood of $x_0$ and $F'$ has compact support, and a closed
convex proper cone $\gamma$ as in {\rm (4a)} such that
$\Phi_{\gamma}(F')\simeq 0$ in a neighborhood of $x_0$.
\item[{\rm (5a)}]
Same condition as {\rm (3a)} with $G\in D^b_{\R-c}(k_X)$.
\enum
\end{lemma}
\begin{proof}
The plan of the proof is as follows:
$$
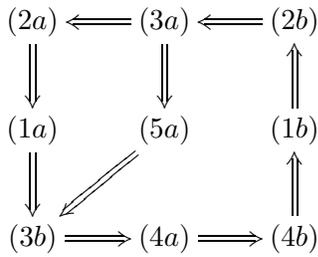
\xymatrix{
(2a)\ar@{=>}[d]&(3a)\ar@{=>}[l]\ar@{=>}[d]&\ar@{=>}[l](2b)&\\
(1a)\ar@{=>}[d]&(5a)\ar@{=>}[dl] &(1b)\ar@{=>}[u]& \\
(3b)\ar@{=>}[r]&(4a)\ar@{=>}[r]&(4b)\ar@{=>}[u]&
}$$

\noindent
(2a) $\Rightarrow$ (1a) follows from
\ $F\simeq Q(\inddlim[i]F_i)\Rightarrow J(F)\simeq \inddlim[i]J(Q(F_i))$.

\medskip
\noindent
(1a) $\Rightarrow$ (3b). 
Let $F\simeq Q(\inddlim[i]F_i)$ and let $i\in I$. There exists $i\to j$ such
that the morphism $F_i\to F_j$ in $D^b(k_X)$ is zero in 
$D^b(k_X;U)$. Hence, there exists a
morphism $F_j\to F'_{ij}$ in $D^b(k_X)$ which is an isomorphism on $U$
and such that the composition 
$F_i\to F_j\to F'_{ij}$ is the zero-morphism in $D^b(k_X)$.
Consider the commutative diagram in which the row on the bottom 
is a distinguished
triangle in $D^b(k_X)$ and $SS(F_{ij})\cap U=\emptyset$:
$$\xymatrix{
            &F_i\ar@{.>}[ld]\ar[d]\ar[rd]|0&                  &\\
F_{ij}\ar[r]&F_j\ar[r]                     &F'_{ij}\ar[r]^{+1}&
}$$
Since the arrow $F_i\to F'_{ij}$ is zero, the dotted arrow may be
completed, making the diagram commutative. 
Hence, we may assume from the beginning that for any $i\in I$ there
exists $i\to j$ such that the morphism $F_i\to F_j$ factorizes as
$F_i\to F_{ij}\to F_j$ with $SS(F_{ij})\cap U=\emptyset$.

We may assume $X$ is affine and $U=W\times \lambda$ where $W$ is open
and relatively compact and $\lambda$ is an open convex cone.
Then $SS(G\tens F_{ij})\cap U=\emptyset$, and
the sheaf $G\tens F_{ij}$ has compact support. Hence, 
$R\sect(X;G\tens F_{ij}) \simeq 0$ which implies 
$H^j\risect(X;G\tens F)\simeq
\inddlim[i]H^j\rsect(X;G\tens F_i) \simeq 0$ for all $j$.
We conclude therefore $\risect(X;G\tens F) \simeq 0$.

\medskip
\noindent
(3b) $\Rightarrow$ (4a). 
Let $F=Q(\inddlim[i]F_i)$, with $F_i\in C^{[a,b]}(\md[k_X])$. 
Set 
$$H_\epsilon=\{x;\langle x-x_0;\xi_0\rangle>-\epsilon\}$$
and let $K\ssubset \pi(U)$ be a compact neighborhood of
$x_0$.
Then there exist an open convex cone $\gamma$ 
and an open neighborhood $W$ of $x_0$ satisfying the following conditions:
\eqn
&&\left\{\ba{l}
W\subset H_\epsilon\cap K,\\
(x+\gamma)\cap H_\epsilon\subset W \mbox{ for all }x\in W,\\
\ol{W}\times\gamma^\circ\subset U\cup T^*_XX.
\ea\right.
\eneqn
Set 
\eqn
&&G_x=k_{(x+\gamma^a)\cap H_\epsilon}, \quad G=\bigoplus_{x\in W}G_x.
\eneqn
Since $\supp(G)\ssubset \pi(U)$ and
$SS(G)\subset \overline W\times{\gamma}^\circ{}^a$, we get
by the hypothesis:
\eqn
\inddlim[i]H^k R\sect(X;G\tens F_i)&\simeq& 0.
\eneqn
Hence, 
\eqn
\inddlim[i](\bigoplus_{x\in W}H^k R\sect(X;G_x\tens F_i)) &\simeq& 0.
\eneqn

Hence one obtains:
\eqn
\left\{ 
\parbox{300pt}{
for any $i\in I$, there exists $i\to j$ such that 
$H^kR\sect(X;G_x\tens F_i)\to H^kR\sect(X;G_x\tens F_j)$ is zero
for any $x\in W$ and any $k\in\Z$.}
\right.
\eneqn
On the other-hand, 
\eqn
&&H^k(\Phi_{\gamma,W}(F_{i}))_x\simeq H^kR\sect(X;G_x\tens F_i).
\eneqn 
Therefore, for any $i\in I$ there exists $i\to j$ such that 
for any $k\in\Z$, the morphism 
$H^k(\Phi_{\gamma,W}(F_{i}))\to H^k(\Phi_{\gamma,W}(F_{j}))$ is the zero
morphism, and this implies
\eqn
&&H^k(\Phi_{\gamma,W}(F))\simeq
\inddlim[i]H^k\Phi_{\gamma,W}(F_{i})
\simeq 0.
\eneqn
This gives the desired result: $\Phi_{\gamma,W}(F)=0$.

\medskip
\noindent
(4a) $\Rightarrow$ (4b) is obvious by taking $F_W$ as $F'$.

\medskip
\noindent
(4b) $\Rightarrow$ (1b).
Let $W$ be an open relatively compact neighborhood of $x_0$
such that $F\vert_W\simeq F'\vert_W$ and $\Phi_\gamma(F')\vert_W\simeq 0$.

Then one has a distinguished triangle:
$$\reeim{q_1}
(k_{\opb{s}(\gamma\setminus\{0\})\cap\opb{q_1}W}\tens \opb{q_2}F')
\to \Phi_\gamma(F')_W\to F'_W\xrightarrow{+1},$$
and hence one obtains
$\reeim{q_1}
(k_{\opb{s}(\gamma\setminus\{0\})\cap\opb{q_1}W}[1]
\tens \opb{q_2}F')\simeq F'_W$.
 Let $F'=Q(\inddlim[i]F_i)$ with $F_i\in C^{[a,b]}(\md[k_X])$, 
and take a finite injective resolution $I$ of 
$k_{\opb{s}(\gamma\setminus\{0\})\cap\opb{q_1}W}[1]$.
Since $I\otimes F_i$ is a finite complex of soft sheaves,
$\reim{q_1}(k_{\opb{s}(\gamma\setminus\{0\})\cap\opb{q_1}W}[1]
\tens \opb{q_2}F_i)$ is represented by
$F'_i:=\eim{q_1}(I\otimes \opb{q_2}F_i)$.
Hence one has
$$\reeim{q_1}
(k_{\opb{s}(\gamma\setminus\{0\})\cap\opb{q_1}W}
\tens \opb{q_2}F')\simeq Q(\inddlim[i]F'_i). $$
Since $SS(F'_i)\cap W\times\Int \gamma^\circ=\emptyset$,
we obtain the desired result.

\medskip
\noindent
(1b) $\Rightarrow$ (2b) is obvious.

\medskip
\noindent
(2b) $\Rightarrow$ (3a). Let $J(F)\simeq \inddlim[i]J(F_i)$. 
If $G\in D^b(k_X)$, we get the isomorphism:
$$\Hom[{D^b(\II[k_X])}](G,F)\simeq\indlim[i]\Hom[D^b(k_X)](G,F_i).$$
We may assume that $X$ is affine and $U=W\times \lambda$ where $W$ is open
and $\lambda$ is an open convex cone. Then the micro-support of
$\rhom(G,F_i)$ is contained in $SS(F_i)+\bar{\lambda}^a$ and
this set does not intersect $X\times\lambda$. Since $\rhom(G,F_i)$
has compact support, $\Hom(G,F_i)$ is zero.

\medskip
\noindent
(3a) $\Rightarrow$ (2a). We may assume that $X$ is affine, $p=(x_0;\xi_0)$ and 
$U=X'\times\Int{\gamma}^\circ$, with $\xi_0\in\Int {\gamma}^\circ$
for a neighborhood $X'$ of $x_0$.
Let $V$ be an open neighborhood of $x_0$ and let 
$W=\{x;\langle x-x_0;\xi_0\rangle>-\epsilon\}$.
Then by taking $V$ and $\epsilon$ small enough,
the sheaf $\Phi_\gamma(H_W)_V$ satisfies
the condition in (3a) for any $H\in D^b(k_X)$.
Let $J(F)=\inddlim[i]J(F_i)$. 
Then 
$\indlim[i]\Hom[D^b(k_X)](G,F_i)\simeq 0$ for any $G=\Phi_\gamma(H_W)_V$. 
Let $i\in I$ and choose $H=F_{i}$.
There exists $i\to j$
such that the composition $(\Phi_\gamma(F_{iW}))_V\to F_i\to F_j$ is zero.
The morphism $(\Phi_\gamma((F_{iW}))_V\to F_i$ is an isomorphism on 
$U':=(V\cap W)\times \Int{\gamma}^\circ$. 
Therefore, $F_i\to F_j$ is zero in $D^b(k_X;U')$.

\medskip
\noindent
(3a) $\Rightarrow$ (5a) is obvious.

%

\medskip
\noindent
(5a) $\Rightarrow$ (3b). (Assuming $F\in D^b_{\IRC}(\II[k_X])$.)
Let (2a-rc) denote the condition (2a) in which one asks moreover that 
$F_i\in D_{\R-c}^{[a,b]}(k_X)$. Define similarly (1a-rc).
Then the same proof of
(3a) $\Rightarrow$ (2a) $\Rightarrow$ (1a) $\Rightarrow$ (3b)
can be applied to show (5a) $\Rightarrow$ (2a-rc) 
$\Rightarrow$ (1a-rc) $\Rightarrow$ (3b).

\end{proof}

\begin{definition}\label{def:ssf}
Let $F\in D^b(\II[k_X])$. The micro-support of $F$, denoted by $SS(F)$, is
the closed conic subset of $T^*X$ whose complementary is the
set of points $p\in T^*X$ such that one of the equivalent conditions in
Lemma~\ref{le:ssf} is satisfied. 
\end{definition}

\begin{proposition}\label{th:sspropert}
\bnum
\item
For $F\in D^b(\II[k_X])$, one has $SS(F)\cap T^*_XX=\supp(F)$.
\item 
Let $F\in D^b(k_X)$. Then $SS(\iota_XF)=SS(F)$.
\item 
Let $F\in D^b(\II[k_X])$. Then $SS(\alpha_XF)\subset SS(F)$.
\item
Let $F_1\to F_2\to F_3\xrightarrow{+1}$ be a distinguished triangle in
$D^b(\II[k_X])$. Then $SS(F_i)\subset SS(F_j)\cup SS(F_k)$ if
$\{i,j,k\}=\{1,2,3\}$.
\enum
\end{proposition}
\begin{proof}
\begin{itemize}
\item[(i)]
$\supp(F)\subset SS(F)$ follows 
for example from (1b) of Lemma~\ref{le:ssf}.
The other inclusion is obvious.

\item[(ii)]
 The inclusion $SS(F)\subset SS(\iota_XF)$ follows from (2a)
since $J(F)$ is $\smash{\inddlim F}$.
The converse inclusion follows from (1b).

\item[(iii)] is obvious, using condition (3b).

\item[(iv)] is obvious by (3b).
\end{itemize}
\end{proof}

\begin{definition}\label{def:limsets}
Let $\Lambda_i, i\in I$ be a family of closed conic subsets of $T^*X$, indexed
by the objects of a small and filtrant category $I$. One sets
\eqn
&&\llim[i]\Lambda_i=\bigcap_{J\subset I}\overline{\bigcup_{j\in J}\Lambda_j}\\
&&\hs{50pt}
\mbox{where $J$ ranges over the family of cofinal subcategories of $I$.}
\eneqn
In other words, $p\in T^*X$ does not belong to 
$\llim[i]\Lambda_i$ if there exists an open neighborhood
$U$ of $p$ and a cofinal subset $J$ of $I$ such that $\Lambda_j\cap
U=\emptyset$ for every $j\in J$.
\end{definition}

It follows immediately from the definition that if
$J(F)\simeq\inddlim[i]J(F_i)$, then 
\eq\label{eq:sslim}
&&SS(F)\subset\llim[i]SS(F_i).
\eneq
It follows from Proposition \ref{th:JtensandJihom} that if $G\in D^b(k_X)$,
one has the inclusions
\eq\label{eq:sslim2}
&&\ba{rcl}
SS(G\tens F)&\subset& \llim[i](SS(G)\widehat{+}SS(F_i)),\\[10pt]
SS(R\ihom(G,F))&\subset& \llim[i](SS(G)^a\widehat{+}SS(F_i)).
\ea
\eneq

\begin{example}\label{exa:sslim}
Let $X=\R^2$ endowed with coordinates $(x,y)$ and denote by 
$(x,y;\xi,\eta)$ the associated coordinates on $T^*X$.
Let
\eqn
Y         &=&\{(x,y);y=0\},\\
U&=&\{(x,y);x^2<y\},\\
Z_\epsilon&=&\{(x,y);x^2<y\leq \epsilon^2 \}.
\eneqn
Set $F_\epsilon=k_{Z_\epsilon}$ and 
$F=k_U\otimes\beta_X(k_{\{0\}})\simeq\inddlim[\epsilon]F_\epsilon$.
Then 
\eqn
SS(k_Y)&=&T^*_YX=\{(x,y;\xi,\eta)\,;\, y=\xi=0\},\\
SS(F_\epsilon)&=&\{(x,y;0,0)\,;\,x^2\le y\leq \epsilon^2\}\\
&&\hs{10pt}\bigcup
\{(x,y;\xi,\eta)\,; \,y=x^2,\left\vert x\right\vert\le\epsilon,
\ \xi=-2x\eta,\ \eta\le0\}\\
&&\hs{30pt}\bigcup
\{(x,y;\xi,\eta)\,;\, y=\epsilon^2,\left\vert x\right\vert\le\epsilon,\ 
\xi=0,\ \eta\le0\}\\
&&\hs{50pt}\bigcup
\{(\pm\epsilon,\epsilon^2;\xi,\eta)\,; \,
0\le\pm\xi\le-2\epsilon\eta,\ \eta\le0\},\\
SS(F)&=&\{(x,y;\xi,\eta); x=y=\xi=0, \eta\leq 0\}.
\eneqn
%
%
On the other-hand, one has
\eqn
&&SS(F)=\llim[\epsilon]SS(F_\epsilon),\\
&&\rhom(k_Y,F)\simeq k_{\{0\}}\left[-2\right],\\
&&\llim[\epsilon](T^*_YX\widehat{+}SS(F_\epsilon))=T^*_{\{0\}}X,\\
&&T^*_YX\widehat{+}SS(F)=\{(x,y;\xi,\eta); x=y=\xi=0\}\\
&&\hs{150pt}\subsetneqq SS(\rhom(k_Y,F)).
\eneqn
Note that $SS(F)$ is not involutive.
\end{example}

Recall that subanalytic isotropic subsets of $T^*X$ are defined in \cite{K-S1}.
Let us say for short that a conic locally closed subset $\Lambda$ of
$T^*X$ is isotropic if $\Lambda$ is contained in a conic locally closed
subanalytic isotropic subset.

\begin{definition}\label{def:wrconst1}
\bnum
\item
We denote by $D^b_{\wRC}(\II[k_X])$ the full triangulated subcategory of 
$D^b_{\IRC}(\II[k_X])$ consisting of objects $F$ such 
that $SS(F)$ is isotropic. We call an object of this category a weakly 
$\R$-constructible ind-sheaf.
\item
Let us denote by $D^b_{\RC}(\II[k_X])$ the full triangulated subcategory of 
$D^b_{\wRC}(\II[k_X])$ consisting of objects $F$ such that 
one has $\rhom(G,F)$ $\in D^b_{\RC}(k_X)$ for any 
$G\in D^b_{\RC}(k_X)$.
We call an object of this category an $\R$-constructible ind-sheaf.
\enum
\end{definition}
Note that the functor $\alpha_X$ induces functors
\eqn
&&\alpha_X\cl D^b_{\wRC}(\II[k_X])\to D^b_{\wRC}(k_X),\\
&&\alpha_X\cl D^b_{\RC}(\II[k_X])\to D^b_{\RC}(k_X).
\eneqn
The last property follows from
$\alpha_X(F)=\rhom(\C_X,F)$.


\begin{conjecture}\label{conj:wrc}
Let $F\in D^b_{\wRC}(\II[k_X])$ and let $G\in D^b_{\wRC}(k_X)$. Then 
$R\ihom(G,F)$ and $G\tens F$ belong to $D^b_{\wRC}(\II[k_X])$.
\end{conjecture}

Example \ref{exa:sslim} shows that the knowledge of $SS(F)$ and $SS(G)$
does not allows us
to estimate the micro-support of $\rhom(F,G)$ by 
the one for sheaves, and that is one reason for the 
definition below.

\begin{definition}\label{def:ssregf}
Let $F\in D^b(\II[k_X])$.
\bnum
\item 
Let $S\subset T^*X$ be a locally closed conic subset and let $p\in T^*X$.
We say that $F$ is regular along $S$ at $p$ if there exist
$F'$ isomorphic to $F$ in a neighborhood of $\pi(p)$, 
an open neighborhood $U$ of $p$ with $S\cap U$ closed in $U$,
a small and
filtrant category $I$ and a functor $I\to D^{[a,b]}(k_X), i\mapsto F_i$ 
such that 
$J(F')\simeq\inddlim[i]J(F_i)$ and $SS(F_i)\cap U\subset S$.
\item
If $U$ is an open subset of $T^*X$ and 
$F$ is regular along $S$ at each $p\in U$, we say that $F$ is
regular along $S$ on $U$.
\item
Let $p\in T^*X$. We say that $F$ is regular at $p$ if $F$ is regular
along $SS(F)$ at $p$. 

If $F$ is regular at each $p\in SS(F)$, 
we say that $F$ is regular.
\item
We denote by $SS_{reg}(F)$ the
conic open subset of $SS(F)$ consisting of points $p$ such that $F$ is
regular at $p$, and we set 
$$  SS_{irr}(F)=SS(F)\setminus SS_{reg}(F). $$
\enum
\end{definition}

Note that $SS_{irr}(F)=SS(F)$ for $F$ in Example \ref{exa:sslim}.
\begin{proposition}\label{th:ssregpropert}
\bnum
\item
Let $F\in D^b(\II[k_X])$. Then $F$ is regular along any locally closed
set $S$ at each $p\notin SS(F)$.
\item
Let $F_1\to F_2\to F_3\xrightarrow{+1}$ be a distinguished triangle in
$D^b(\II[k_X])$. 
If $F_j$ and $F_k$ are regular along $S$, so is $F_i$ 
 for $i,j,k\in\{1,2,3\}, j\neq k$.
\item 
Let $F\in D^b(k_X)$. Then $\iota_XF$ is regular.
\enum
\end{proposition}
\begin{proof}
(i) and (iii) are obvious and 
the proof of (ii) is similar to that of Proposition~\ref{th:sspropert}~(iv).
\end{proof}


It is possible to localize the category $D^b(\II[k_X])$ with respect
to the micro-support, exactly as for usual sheaves.

Let $V$ be a subset of $T^*X$ and let $\Omega=T^*X\setminus V$. 
We shall denote by $D^b_V(k_X)$ the full triangulated 
subcategory of $D^b(k_X)$ consisting
of objects $F$ such that $SS(F)\subset V$, and by 
$D^b(k_X;\Omega)$ the localization of 
$D^b(k_X)$ by $D^b_V(k_X)$. 

Similarly, we denote
by $D^b_V(\II[k_X])$ the full triangulated 
subcategory of $D^b(\II[k_X])$ consisting
of objects $F$ such that $SS(F)\subset V$.

\begin{definition}
One sets 
$$D^b(\II[k_X;\Omega])=D^b(\II[k_X])/D^b_V(\II[k_X]),$$
the localization of $D^b(\II[k_X])$ by $D^b_V(\II[k_X])$.
\end{definition}
Let $F_1$ and $F_2$ are two objects
of  $D^b(\II[k_X])$ whose images in  $D^b(\II[k_X;\Omega])$ are
isomorphic. There exist a third object
$F_3\in D^b(\II[k_X;\Omega])$ and
distinguished triangles in $D^b(\II[k_X])$:
 $F_i\to F_3\to G_i\xrightarrow{+1}$ ($i=1,2$) 
such that $SS(G_i)\cap \Omega=\emptyset$. 
It follows that $SS(F_1)\cap\Omega=SS(F_3)\cap\Omega=SS(F_2)\cap\Omega$.

Therefore if $F\in D^b(\II[k_X;\Omega])$,
the subsets $SS(F)$
and $SS_{irr}(F)$ of $\Omega$
are well-defined. 

\section{Invariance by contact transformations}

It is possible to 
define contact transformations on ind-sheaves. We shall follow the notations
in \cite{K-S1} Chapter~VII. 

We denote by $p_1$ and $p_2$ the first and second projections defined
on $T^*(X\times Y)\simeq T^*X\times T^*Y$, and we denote by $p_2^a$
the composition of $p_2$ with the antipodal map on $T^*Y$. 

We denote by $r\cl X\times Y\to Y\times X$ the canonical map and we keep
the same notation to denote its inverse.

By a kernel $K$ on $X\times Y$ 
we mean an object of $D^b(k_{X\times Y})$. To a kernel
$K$ one associates the kernel on $Y\times X$ 
\eqn
K^*:=\oim{r}\rhom(K,\omega_{X\times Y/Y}).
\eneqn
One defines the functor
\eq\label{eq:phik1}
\Phi_K\cl  D^b(k_Y)&\to& D^b(k_X)\\
G&\mapsto& R\eim{q_1}(K\tens\opb{q_2}G).\nonumber
\eneq 

Consider another manifold $Z$ and a kernel $L$ on $Y\times Z$.
One defines the projection $q_{12}$ from $X\times
Y\times Z$ to $X\times Y$, and similarly with  $q_{23}$,  $q_{13}$.

One sets
\eq\label{eq:circkernels}
K\circ L=R\eim{q_{13}}(\opb{q_{12}}K\tens \opb{q_{23}}L).
\eneq
Choosing $Z=\{\rmpt\}$, one has 
$\Phi_K(G)=K\circ G$ for $G\in D^b(k_Y)$.

Let $\Omega_X$ and $\Omega_Y$ be two conic open subsets of $T^*X$ and
$T^*Y$, respectively.
One denotes  by $N(\Omega_X,\Omega_Y)$ the full subcategory of 
$D^b(k_{X\times Y};\Omega_X\times T^*Y)$ of objects $K$ satisfying;
\eq\label{eq:catnomega}
\left\{ \begin{array}{l}
SS(K)\cap (\Omega_X\times T^*Y)\subset \Omega_X\times\Omega_Y^a,\\[5pt]
p_1\cl SS(K)\cap(\Omega_X\times T^*Y)\to \Omega_X \mbox{ is proper.}
\end{array}\right.
\eneq

Let us recall some results of loc. cit.
\bnum
\item
Let $K\in N(\Omega_X,\Omega_Y)$. Then the functor $\Phi_K$ induces a
well-defined functor:
$\Phi_K^\mu\cl  D^b(k_Y;\Omega_Y)\to D^b(k_X;\Omega_X)$.
\item
Let $L\in N(\Omega_Y,\Omega_Z)$. 
Then $K\circ L\in N(\Omega_X,\Omega_Z)$.
Moreover, the two functors 
$\Phi_{K\circ L}^\mu$ and $\Phi_K^\mu \circ \Phi_L^\mu$
from $D^b(k_Z;\Omega_Z)$ to $D^b(k_X;\Omega_X)$ are isomorphic.
\enum

We construct the functor analogous to the functor $\Phi_K$ for ind-sheaves 
by defining

\eq\label{eq:phik2}
\tPhi_K\cl  D^b(\II[k_Y])&\to& D^b(\II[k_X])\\
G&\mapsto& R\eeim{q_1}(K\tens\opb{q_2}G).\nonumber
\eneq 

Applying Theorem \ref{th:indcab2}, we get:
\begin{lemma}\label{le:indker1}Let $G\in D^b(\II[k_Y])$ and assume
  that 
$J(G)\simeq \inddlim[i]J(G_i)$, with $I$ small and filtrant and $G_i\in
D^b(k_Y)$. Then 
$J(\tPhi_K(G))\simeq \inddlim[i]J(\Phi_K(G_i))$.
\end{lemma}

Now assume that $\dim X=\dim Y$ and 
that there exists a smooth conic Lagrangian submanifold 
$\Lambda\subset \Omega_X\times\Omega_Y^a$ such that 
$p_1\cl \Lambda\to \Omega_X$ and $p_2^a\cl \Lambda\to \Omega_Y$ are
isomorphisms. In other words, $\Lambda$ is the graph of a homogeneous 
symplectic isomorphism $\chi\cl  \Omega_Y\isoto\Omega_X$.

Let $K$ be a kernel satisfying the assumptions of Theorem~7.2.1 of
loc. cit., that is:
\eq\label{eq:kernelqct}
\left\{ \begin{array}{l}
\mbox{$K$ is cohomologically constructible,}\\[5pt]
(\opb{p_1}(\Omega_X)\cup\opb{{p_2^a}}(\Omega_Y))\cap SS(K)
\subset \Lambda,\\[5pt]
k_\Lambda\isoto\mu hom(K,K)\quad \mbox{on $\Omega_X\times\Omega_Y^a$.}
\end{array}\right.
\eneq

\begin{theorem}\label{th:indker1}
Assume \eqref{eq:kernelqct}.
\bnum
\item
The functor $\tPhi_K$ induces a
well-defined functor:
$\tPhi_K^\mu\cl D^b(\II[k_Y;\Omega_Y])\to D^b(\II[k_X;\Omega_X])$.
Similarly, the functor $\tPhi_{K^*}$ induces a
well-defined functor:
$\tPhi_{K^*}^\mu\cl D^b(\II[k_X;\Omega_X])\to D^b(\II[k_Y;\Omega_Y])$.
\item
The functor
$\tPhi_K^\mu\cl D^b(\II[k_Y;\Omega_Y])\to D^b(\II[k_X;\Omega_X])$
and the functor
$\tPhi_{K^*}^\mu\cl  D^b(\II[k_X;\Omega_X])\to D^b(\II[k_Y;\Omega_Y])$
are equivalences of
categories inverse one to each other.
\item
If $G\in D^b(\II[k_Y])$, then 
$SS(\tPhi_K(G))\cap\Omega_X=\chi(SS(G)\cap \Omega_Y)$.
\item
If $G$ is regular at $p\in \Omega_Y$, then $\tPhi_K(G)$ is regular 
at $\chi(p)\in \Omega_X$. In other words,
$SS_{irr}(\Phi_K(G))\cap\Omega_X=\chi(SS_{irr}(G)\cap \Omega_Y)$.
\enum
\end{theorem}
\begin{proof}
(i) Let $G\in D^b(\II[k_Y])$ and assume that 
$SS(G)\cap\Omega_Y=\emptyset$. Let us prove that 
$SS(\tPhi_K(G))\cap\Omega_X=\emptyset$. 
Let $p_X\in\Omega_X$ and let $p_Y=\opb{\chi}(p_X)$. There
exist an open neighborhood $U_Y$ of $p_Y$ in $\Omega_Y$ and an
inductive system such that $J(G)\simeq \inddlim[i\in I]J(G_i)$, and for
any $i\in I$ there exists $i\to j$ such that the morphism 
$G_i\to G_j$ is zero in $D^b(k_Y;U_Y)$. 
Applying Lemma \ref{le:indker1} we find that 
$J(\tPhi_K(G))\simeq \inddlim[i]J(\Phi_K(G_i))$. Since the morphism 
$\Phi_K(G_i)\to \Phi_K(G_j)$ is zero in $D^b(k_X;U_X)$, the result follows.

\noindent
(ii) One has the isomorphism $K\circ K^*\simeq k_{\Delta_X}$ 
in $N(\Omega_X,\Omega_X)$ and the isomorphism  
$K^*\circ K\simeq k_{\Delta_Y}$ in $N(\Omega_Y,\Omega_Y)$.
Hence, it is enough to remark that 
\eq\label{eq:kcirckmu}
&&\tPhi_K^\mu\circ\tPhi_{K^*}^\mu\simeq \tPhi_{K\circ K^*}^\mu,
\eneq
which follows from the fact that the two functors
$\tPhi_K\circ\tPhi_{K^*}$ and $\tPhi_{K\circ K^*}$, from 
$D^b(\II[k_X])$ to $D^b(\II[k_X])$ are isomorphic.

\noindent
(iii) For an open subset $U_Y\subset\Omega_Y$, set $U_X=\chi(U_Y)$. Then 
$K\in N(U_X,U_Y)$ and $K$ satisfies \eqref{eq:kernelqct} with
$\Omega$ replaced with $U$. Let $G\in D^b(\II[k_Y])$ 
with $SS(G)=\emptyset$ in a neighborhood of $p_Y\in\Omega_Y$. 
By the proof of (i), 
$SS(\tPhi_K(G))=\emptyset$ in a neighborhood of $\chi(p_Y)$.

\noindent
(iv) The proof is similar to that of (iii).
\end{proof}

\section{Ind-sheaves and $\shd$-modules}

Let now $X$ be a complex manifold and 
let $\shm$ be a coherent $\shd_X$-module.
We set for short
\eqn
Sol(\shm)&=& \rhom[\shd_X](\shm,\sho_X),\\
Sol^t(\shm)&=&R\ihom[\beta_X\shd_X](\beta_X\shm,\sho_X^t).
\eneqn

\begin{theorem}\label{th:ssdmod}
One has 
\eqn
SS(Sol^t(\shm))=\car(\shm).
\eneqn
\end{theorem}
\begin{proof}
(i) The inclusion $\car(\shm)\subset
SS(Sol^t(\shm))$ follows from  
$$SS(Sol(\shm))=\car(\shm),\quad\alpha_X(Sol^t(\shm)) \simeq Sol(\shm).$$
and Proposition~\ref{th:sspropert}~(ii).

\medskip
\noindent
(ii) Let us prove the converse inclusion using condition (5a) of 
 Lemma~\ref{le:ssf}. 
Assume that $G\in D^b_{\R-c}(\C_X)$ satisfy
$SS(G)\cap\car(\shm)\subset T^*_XX$.
One has the morphisms
\eqn
R\hom(G,R\ihom[\beta_X\shd_X](\beta_X\shm,\sho_X^t))&\simeq&
R\hom[\shd_X](\shm,\thom(G,\sho_X))\\
&\to& R\hom[\shd_X](\shm,R\hom(G,\sho_X)).
\eneqn
It follows from  \cite[Corollary~4.2.5]{An} that 
the second morphism is an isomorphism.
Hence the result follows from
$SS(Sol(\shm))=\car(\shm)$ and  Lemma~\ref{le:ssf} (5a).
\end{proof}

The following conjecture is a consequence of Conjecture \ref{conj:wrc}.

\begin{conjecture}
If $\shm$ is a holonomic $\shd_X$-module, then
$Sol^t(\shm)$ belongs to $D^b_{\RC}(\II[\C_X])$.
\end{conjecture}


\begin{theorem}\label{th:regdmod}
If $\shm$  is a regular holonomic $\shd_X$-module, then
$Sol^t(\shm)\to Sol(\shm)$ is an isomorphism.
\end{theorem}
\begin{proof}
This is a reformulation
of a result of \cite{K1} which asserts that for any $G\in
D^b_{\R-c}(\C_X)$, the natural morphism 
\eqn
\rhom[\shd_X](\shm,\thom(G,\sho_X))\to
\rhom[\shd_X](\shm,\rhom(G,\sho_X))
\eneqn
is an isomorphism.
%
\end{proof}

We conjecture the following statement in which
``only if '' part is a consequence of the theorem above.
\begin{conjecture}
Let $\shm$  be a holonomic $\shd_X$-module. Then $\shm$ is regular
holonomic if
and only if $Sol^t(\shm)$ is regular.
\end{conjecture}

\section{An example}

In this section $X=\C$ endowed with the holomorphic coordinate $z$, 
 and we shall study the ind-sheaf of temperate holomorphic
solutions of the $\shd_X$-module $\shm:=\shd_X\exp(1/z)
=\shd_X/\shd_X(z^2\partial_z+1)$.
We set for short
\eqn
\ba{rll}
\shs^t{:=}&H^0(Sol^t(\shm))&\simeq\ihom[\beta_X\shd_X](\beta_X\shm,\sho_X^t),
\\[5pt]
\shs {:=}&H^0(Sol(\shm))&\simeq\hom[\shd_X](\shm,\sho_X).
\ea
\eneqn

Notice first that $\sho_X^t$ is concentrated in degree $0$ (since
$\dim X=1$), and it is a sub-ind-sheaf of $\sho_X$. It follows that the
morphism $\shs^t\to\shs$ is a monomorphism.

Moreover, 
\eqn
&&\shs\simeq\C_{X,X\setminus \{0\}}\cdot\exp(1/z).
\eneqn

\begin{lemma}
Let $V\subset X$ be a connected open subset. Then 
$\sect(V;\shs^t)\neq 0$ if and only if $V\subset X\setminus \{0\}$ and 
$\exp(1/z)\vert_V$ is tempered.
\end{lemma}
\begin{proof}
The space $\sect(V;\shs)$ has dimension one and is generated by the
function $\exp(1/z)$. Hence, the subspace 
$\sect(V;\shs^t)\simeq \sect(V;\shs)\cap\sect(V;{\sho}^t)$ is not zero
if and only if $\exp(1/z)\in \sect(V;\sho_X^t)$, that is, if and only
if $\exp(1/z)\vert_V$ is tempered.
\end{proof}
Let us set $z=x+iy$. 

\begin{lemma}\label{le:expbnd}
Let $W$ be an open subanalytic subset of $\BBP^1(\C)$ with
$\infty\notin W$. Assume that 
there exist positive constants $C$ and $A$ such that 
\eq\label{eq:exp1z}
&&\exp(x)\leq  C(1+x^2+y^2)^N\mbox{ on }W.
\eneq
Then there exists a constant $B$ such that 
$x\leq B$ on $W$.
\end{lemma}
\begin{proof}
If $x$ is not bounded on $W$, then there
exists a real analytic curve $\gamma\cl [0,\epsilon[\to\BBP^1(\C)$ such
that $\Real \gamma(0)=\infty$ and $\gamma(t)\in W \mbox{ for }t>0$.
Writing $\gamma(t)=(x(t), y(t))$, 
one has
$$y(t)=cx(t)^q+\mathrm{O}(x(t)^{q-\epsilon}).$$
for some $q\in \Q$, $c\in\R$ and $\epsilon>0$.
Then \eqref{eq:exp1z} implies that $\exp(x)$ has a polynomial growth
when $x\to\infty$, which is a contradiction.
\end{proof}

Let $\bar B_\epsilon$ denote the closed ball with center $(\epsilon,0)$ 
and radius $\epsilon$ and set 
$U_\epsilon=X\setminus \bar B_\epsilon$.

\begin{proposition}
One has the isomorphism
\eq\label{eq:solt=indlim}
&&
\inddlim[\epsilon>0]\C_{XU_\epsilon}\isoto
\ihom[\beta_X\shd_X](\beta_X\shm,\sho_X^t).
\eneq
\end{proposition}
\begin{proof}
It follows from Lemma \ref{le:expbnd} that 
$\exp(1/z)$ is temperate (in a neighborhood of $0$) 
 on an open subanalytic subset $V\subset X\setminus\{0\}$ if and only
if ${\rm Re}(1/z)$ is bounded on $V$, that is, if and only if 
$V\subset U_\epsilon$ for some $\epsilon>0$.

Let $V$ be a connected relatively compact subanalytic open subset of
$X\setminus\{0\}$. Then a morphism 
$\C_{V}\to\C_{X\setminus\{0\}}\cdot\exp(1/z)$ factorizes through a morphism
 $\C_{V}\to \shs^t$ if and only if
it factorizes through $\C_{U_\epsilon}$. 
Hence
we get the isomorphism \eqref{eq:solt=indlim}
by Theorem \ref{th:suba}.
\end{proof}

\begin{remark}
In fact one can show
$$H^1(Sol^t(\shm))\isoto H^1(Sol(\shm))\simeq\C_{0}.$$
The isomorphism
$H^1(Sol(\shm))=\sho_X/(z^2\partial_z+1)\sho_X\isoto \C_{0}$
is given by
$$(\sho_X)_0\ni v(z)\mapsto \oint v(z)z^{-2}\exp(-1/z)\,dz.$$
Note that $\phi(z){:=}z^{-2}\exp(-1/z)$
is a solution to the adjoint equation
$$(-\partial_zz^2+1)\phi(z)=0.$$
The distinguished triangle
$$\shs^t\to Sol^t(\shm)\to H^1(Sol^t(\shm))[-1]\xrightarrow{+1}$$
gives a non-zero element of
$Ext^2(\C_{0},\shs^t)\isoto Ext^2(\C_{0},\C_X)\simeq\C$.
\end{remark}

\small
Masaki Kashiwara\\
Research Institute for Mathematical Sciences\\
Kyoto University\\
Kyoto 606-8502\\
Japan\\
masaki@kurims.kyoto-u.ac.jp\\
\vspace{5mm}

Pierre Schapira\\
Universit{\'e} Pierre et Marie Curie\\
Case 82, Analyse Alg{\'e}brique, UMR7586\\
4, Pl. Jussieu 75252 Paris Cedex 05\\
France\\
schapira@math.jussieu.fr\\
http://www.math.jussieu.fr/{\~{}} schapira/\\

\end{document}